\documentclass[12pt]{amsart}
\usepackage[T1]{fontenc}
\usepackage{amssymb}
\newcommand{\cl}{\operatorname{cl}}

\newcommand{\Cee }{\mathcal C}
\newcommand{\Wee }{\mathcal W}

\newcommand{\Fee }{\mathcal F}
\newcommand{\Mee }{\mathcal M}
\newcommand{\Uee }{\mathcal U}
\newcommand{\Hee }{\mathcal H}

\newcommand{\Ree }{\mathbb R}
\newcommand{\add}{\operatorname{\textit{add\,}}}
\newcommand{\non}{\operatorname{\textit{non\,}}}
\newcommand{\cov}{\operatorname{\textit{cov\,}}}
\newcommand{\cof}{\operatorname{\textit{cof\,}}}
\newcommand{\nwd}{\operatorname{\textit{nwd}}}
\renewcommand{\int}{\operatorname{Int}}

\newtheorem{theorem}{Theorem}
\newtheorem{corollary}[theorem]{Corollary}
\newtheorem{lemma}[theorem]{Lemma}
\newtheorem*{problem}{Question}

\author{Andrzej Kucharski}
\address{Andrzej Kucharski \\
 Institute of Mathematics, University of
Silesia \\
 Bankowa 14, 40-007 Katowice}
\email{akuchar@ux2.math.us.edu.pl}
\author{Szymon Plewik}
\address{Szymon Plewik\\Institute of Mathematics,
University of Silesia , Ban\-ko\-wa 14, 40-007 Katowice}
\email{plewik@ux2.math.us.edu.pl}

\title{Cardinal invariants for $\Cee$-cross  topologies  }

\begin{document}
              \begin{abstract}
$\Cee$-cross topologies are introduced. Modifications of the Kuratowski-Ulam Theorem  are considered.  Cardinal invariants $\add$, $\cof$, $\cov$ and  $\non$ with respect to meager or nowhere dense subsets are compared. Remarks on invariants $\cof(\nwd_Y)$  are mentioned for  dense  subspaces $Y\subseteq X$. 
               \end{abstract}
               \maketitle
               \section{Introduction}
Let  $X\times Y=\{(x,y): x\in X \mbox{ and } y\in Y\}$ denote the Cartesian product of sets $X$ and $Y$. For a subset  $G \subseteq X \times Y$  let $G_x$ denote a vertical section 
$\{y\in Y: (x,y)\in G\}$ and let $G^y$ denote a horizontal section $\{x\in X: (x,y)\in G\}$.  Let $X$ be a set and $\Fee$  be a family of subsets of $X$. Assume  that $X = \bigcup \Fee$. Consider following cardinal invariants assign to  $X$ and $\Fee$. \\ $ \bullet $ The minimal cardinality of  a subset $Y\subseteq X$ with $Y \notin \Fee$  is denoted by
    $\non(\Fee)$. \\  $ \bullet $ 
 The minimal cardinality of  a subfamily $\Wee\subseteq \Fee$ such that $\bigcup \Wee \notin \Fee$ is denoted by $\add(\Fee)$.\\  $\bullet $ 
  The minimal cardinality of  a subfamily $\Wee\subseteq \Fee$ such that $\bigcup \Wee = X$  is denoted by  $\cov( \Fee )$. \\ $ \bullet $ 
The minimal cardinality of  a subfamily $\Wee\subseteq \Fee$ such that for any  $Y \in \Fee$ there exists $Z\in \Wee$ with $Y\subseteq Z$, which is denoted by  $\cof(\Fee)$. \\
Recall   a few well know (a folklore) facts.\\$ \bullet $ 	 $\cov(\Fee)\leq \cof (\Fee)$. \\
  $ \bullet $ 	 If $X\notin  \Fee$, then $ \non(\Fee)\leq \cof (\Fee).$ \\
 $ \bullet $ 	  If
$X\notin\Fee$, then $\add(\Fee)\leq \cov(\Fee)$. \\
 $ \bullet $ 	If 
$\{ \{x\}: x \in X \} \subseteq \Fee$, then $  \add ( \Fee )\leq \non( \Fee )$.\\ 
Such  cardinal invariants have been used by many authors. Usually, these invariants are considered for  subsets or subfamilies of the reals $\Ree$. Sometimes for other topological spaces. For instance J. Kraszewski \cite{kr} studied inequalities between such cardinal invariants for Cantor cubes.   Z. Piotrowski and A. Szyma\'nski \cite{ps} considered cardinal functions $\add$, $\cov$ and $\non$ for  arbitrary topological spaces. Nevertheless, compare the last survey article by A. Blass \cite{bla} which contains a huge bibliography.

Let $X$ be a topological space with a topology $\mu$. The ideal of all nowhere dense subset of $X$ is denoted by $\nwd_X$ or by $\nwd_\mu$; but   $\Mee_X$ or $\Mee_\mu$ denote  the $\sigma-$ideal of all meager   subsets of $X$. 
A family $\Fee$ of non-empty subsets of  $X$ is a $\pi$\textit{-network} if each  non-empty open subset of
$X$ contains a member of $\Fee$. Any $\pi$-network is called $\pi$\textit{-base} whenever it contains open sets, only.
The minimal cardinality of  a $\pi$-base for $X$ is denoted by  $\pi (X)$. In this note it is  assumed that any cardinal invariant is infinite.   A topological space which satisfies the Baire category theorem is called Baire space.
If a $T_1$-space $X$ is  dense in itself, then   $\bigcup\nwd_X=X$ and $X \notin \nwd_X$ and $\{ \{x\}: x \in X \} \subseteq \nwd_X$. Additionally, if  $X$ is a Baire space,  then   $X\notin \Mee_X$ and $\cov(\nwd_X) \geq \omega_1$. 
                   
For   topological spaces $X$ and $Y$ let $\gamma$ be a family of all subsets of $X \times Y$  such that for each $G \in \gamma$  all vertical sections $G_x$ are open  in $Y$ and  all horizontal sections $G^y$ are open in $X$. The family $\gamma$ is a topology on $X \times Y$. It is usually called the cross topology. A function $f: X\times Y \to Z$ is \textit{separately continuous} whenever  all maps (which   make one of the two variables constant) $y\to f(z,y)$ and $x\to f(x,z)$ are  continuous. In fact, the cross topology   $\gamma$ is the weak topology on $X\times Y$ generated by the family of all separately continuous  functions into any topological space. 
  The product topology on $X\times Y$  is denoted by  $\tau$. The topology $\tau$ is generated by the family of rectangles $\{U\times V: U \subseteq X \mbox{ is open and } V \subseteq Y \mbox{ is open} \}$. The topology of separate continuity  on  $X\times Y$ is  denoted by   $\sigma$. The topology $\sigma$ is the weak topology on $X\times Y$ generated by the family of all separately continuous real-valued functions. The topology $\gamma$ is finer than the product topology $\tau$ and the topology of separate continuity $\sigma$. Topologies $\tau$, $\sigma$ and $\gamma$ have been compared a few time in the literature. The first systematic study was done by C. J. Knight, W. Moran and J. S. Pym \cite{kmp} and \cite{kmp2}, see also J. E. Hart and K. Kunen \cite{hk} or  M. Henriksen and R. G. Woods \cite{hw} for some comments and exhaustive references.

We consider some  modifications  of the cross topology $\gamma$. 
It is introduced the family of $\Cee$-cross topologies on $X\times Y$. Under some additional   assumptions any $\Cee$-cross topology fulfills the Kuratowski-Ulam Theorem, see Theorem \ref{EEE}. A proof that each open and dense set with respect to $\tau$ is  open and dense  with respect to $\Cee$-cross topology needs some  other  assumptions about  $\Cee$, see Corollary \ref{AAA}. Special cases when the product topology  of two Baire space is a Baire space are generalized by Theorem \ref{BBB}.   In \cite{hw} was stated conditions under which each nowhere dense set in $\tau$ is nowhere dense in $\gamma$. Theorem \ref{HHH} and Corollary \ref{FFF} establish conditions  with $\Cee$-cross topologies on $\Ree \times \Ree$ for similar results. Cardinal invariants for various $\Cee$-cross topologies are calculated in theorems \ref{GGG}, \ref{JJJ} and \ref{KKK}. In Lemma \ref{lm2}, Corollary \ref{nwd} and Theorem \ref{4XX} are improved some results concerning cofinality of nowhere dense or meager sets in separable metric spaces, compare \cite{bhh}. 
                    \section{$\mathcal C$-cross topologies} 
Let $X$ and $Y$ be topological spaces. Consider  a family $\mathcal C$ of subsets of $X \times Y$ which is closed under finite intersections and such that for each $G \in \mathcal C$ all vertical sections  $G_x$ have non-empty interior in $Y$, and all horizontal sections  $G^y$ have non-empty interior in $X$. The topology generated by $\mathcal C$ is   called $\mathcal C$\textit{-cross topology}. Obviously, $\mathcal C$ is a base for  the $\mathcal C$-cross topology and topologies $\tau$, $\sigma$ and $\gamma$ are $\mathcal C$-cross topologies. 

Let $\mathbb S$ be the Sorgenfrey line. By $\mathbb S^n$ and $\Ree^n$  denote  products  of $n$-copies of  $\mathbb S$  and the reals $\Ree$ 
 with  product topologies, respectively. One can  check that $  \nwd_{\mathbb S^n} = \nwd_{\Ree^n}$.  
 Indeed, the family of all products of  $n$-intervals with rational endpoints is a $\pi$-base for  $\mathbb S^n$ and $\Ree^n$. So, both spaces have the same nowhere dense subsets.  
 Recall that, see A. Todd \cite{to} or compare  \cite{hw},   topologies $\sigma$ and $\tau$ on a set $X$ are  $\Pi$\textit{-related} if $\tau$ contains a $\pi$-network for the topology $\sigma$ on $X$ and $\sigma$ contains a $\pi$-network for the topology $\tau$ on $X$.
Note  that, different $\mathcal C$-cross topologies on $X\times Y$ have not to be $\Pi$-related, see \cite{hw} Theorem 3.6  for suitable counter-examples. The next lemma is a small modification of Proposition 3.2 in \cite{hw}.
\begin{lemma} \label{XXX} If topologies $\sigma$ and $\tau$ on a set $X$ are $\Pi$-related, then $\sigma$-dense subspaces and $\tau$-dense subspaces  are the same. Also,   $\sigma$-nowhere dense subsets of $X$ and $\tau$-nowhere dense subsets of $X$ are the same. \qed
\end{lemma} 
 Lemma \ref{XXX} generalizes the fact that   spaces  $\mathbb S^n$ and $\Ree^n$ have the same dense subspaces, and the same nowhere dense sets. In \cite{hw}, it is applied to the topology of separate continuity. The plane $\Ree \times \Ree$  has the same dense subspaces and the same nowhere dense sets for the product topology and for the topology of separate continuity, respectively. This is not true for the cross topology, see \cite{hw} Theorem 3.6.a.

The proof of the following theorem requires that $\Cee$ contains a $\pi$-network for the product topology  $\tau$ and it needs the inequality $\pi (Y) < \add (\Mee_X)$. 
          \begin{theorem}\label{EEE} Let $X\times Y$ be equipped with a $\Cee$-cross topology such that the family $\Cee$ contains a $\pi$-network for the product topology  $\tau$   and let $\pi (Y) < \add (\Mee_X)$. If  a set $E\subseteq X\times Y$ is open and dense with respect to the $\mathcal C$-cross topology, then there exists a meager subset $P\subseteq X$ such that   any section $E_x= \{ y\in Y: (x,y) \in E \}$ is dense in $Y,$ for  each $x \in X \setminus P$.
         \end{theorem}
          \begin{proof} Let $\Uee$ be a $\pi$-base for $Y$ of the cardinality $\pi (Y)$. For each $V \in \Uee$ and any non-empty and open $W\subseteq X$ the rectangle $W\times V$ has non-empty interior with respect to  the $\mathcal C$-cross topology, since the family $\Cee$ contains a $\pi$-network for the product topology.  The intersection  $E \cap (W\times V)$ contains a non-empty    $\mathcal C$-cross open  subset. By the definition of a $\mathcal C$-open set there exists $q\in V$ such that $\emptyset \not= \int_X (E \cap(W \times V))^q \subseteq W.$ For abbreviation let  $A^V= \{x\in X: (\{x \} \times V) \cap E \not= \emptyset \} .$ 
Since $W$ is an arbitrary open subset  of $X$ and $ \int_X (E \cap(W \times V))^q \subseteq  A^V ,$ then one concludes that each set $A^V \subseteq X $ contains an open dense subset. Put $P= X \setminus \bigcap \{ A^V: V\in \Uee \}$.  The set $P\subseteq X $ is meager, since $\pi (Y) < \add (\Mee_X)$. For each point $x \in \bigcap \{ A^V: V\in \Uee \}$ the set $ \int_Y E_x \subseteq Y$ has to be dense.
           \end{proof} 
Theorem \ref{EEE} is  a modification of the Kuratowski-Ulam Theorem, compare \cite{ku}, \cite{kur}, \cite{ox} or \cite{fnr}.  One can call it the Kuratowski-Ulam Theorem, too. The next corollary follows from the Kuratowski-Ulam Theorem which is applied to the product topology $\tau$. Its proof  needs that $\Cee$ contains a $\pi$-base instead of a $\pi$-network.
              \begin{corollary}\label{AAA} Let $X\times Y$ be equipped with a $\Cee$-cross topology such that the family $\Cee$ contains a $\pi$-base for the product topology  $\tau$   and let $\pi (Y) < \add (\Mee_X)$ and assume that any  non-empty and open subset of $X$ is not meager.  If  a set $E\subseteq X\times Y$ is open and dense with respect to the product topology, then $E$ contains a subset $G$ which is dense and open  with respect to any $\mathcal C$-cross topology. 
           \end{corollary}
            \begin{proof}  Let $\Uee \subseteq \mathcal C$ be a $\pi$-base for the product topology $\tau$. The union $$G= \bigcup \{ V\in \Uee: V\subseteq E\}$$ is a $\mathcal C$-cross open  set and open dense with respect to the product topology. Apply  the Kuratowski-Ulam Theorem for the product topology $\tau$. There exists a meager subset $P\subseteq X$ such that   the  section $G_x\subseteq Y $ is open and dense, for any  $x \in X \setminus P$. Suppose that a non-empty $\mathcal C$-cross open set $V$ is disjoint with $G$. For any $(p,q) \in V$ we get  $\int_X V^q \not= \emptyset$. So, for each $x\in \int_X V^q$ the set $G_x \subseteq Y$ is not dense.  But the non-meager  set $\int_X V^q$ can  not be contained in the meager set $P$, a contradiction. 
          \end{proof}
Recall that, a space is called \textit{quasiregular} if each non-empty open set contains the closure of some non-empty open set.
Let a space  $Y$  be quasiregular and  $\tau \subseteq \mathcal C$ and $\mu$ be the $\Cee$-cross topology. Suppose that for any   set $G\in \mu$ all  vertical sections $G_x \subseteq Y$ are open. Under such assumptions Corollary \ref{AAA} would be deduced from the proof of Lemma 3.4 (a) in \cite{hw}. One should  adopt the proof, since the lemma concerns Tychonoff spaces and takes the topology $\gamma$ instead of $\mu$. In fact,  $\emptyset \not= V \in \mu$ implies $ \int_\tau \cl_\mu V \not= \emptyset$.  But  $ \mu$ is finer than $ \tau$, hence   $V\cap E =\emptyset$ imply  $\cl_\mu V \cap E = \emptyset$, for each set $E \in \tau$. If additionally $E$ 
 is dense with respect to $\tau$, then it should  be $V =\emptyset$. So,  $E$ has to be dense with respect to $\mu$, too. 
                     \section{$\mathcal C$-meager sets}
Let us examine the Baire category theorem with respect to $\mathcal C$-cross topologies. The next theorem is related to results  which were  obtained  by A. Kucia \cite{kuc} or D. Gauld, S. Greenwood and Z. Piotrowski \cite{gp}. Our proof  is a small improvement of Theorem \ref{EEE}.  
\begin{theorem}\label{BBB} Let $X\times Y$ be equipped with a $\Cee$-cross topology such that the family $\Cee$ contains a $\pi$-network for the product topology  $\tau$ and let $\pi (Y) < \cov (\Mee_X)$.   If a family  $\{E_\alpha \subseteq X\times Y:\alpha<\lambda \}$ consists of sets which are  open and dense with respect to the  $\mathcal C$-cross topology, then the intersection $\bigcap \{E_\alpha : \alpha < \lambda \} $ is non-empty for any cardinal number $\lambda <\min\{\cov (\Mee_X), \cov (\Mee_Y)\}$. 
                 \end{theorem}
               \begin{proof}  Assume that $\Uee$ is a $\pi$-base for $Y$ of the cardinality $\pi (Y)$. If $V \in \Uee$, then for every non-empty and open $W\subseteq X$ the rectangle $W\times V$ has non-empty interior with respect to  the $\mathcal C$-cross topology. 
Similarly like in the proof of Theorem \ref{EEE} one concludes that each set $$A_\alpha^V=\{x\in X: (\{x \} \times V) \cap E_\alpha \not= \emptyset \} $$ contains a  dense open subset of $X$.  There exists a point $$x\in   \bigcap \{ A_\alpha^V: V\in \Uee \mbox{ and } \alpha< \lambda \},$$ since $\pi(Y) < \cov (\Mee_X)$ and $ \lambda < \cov (\Mee_X)$. Any vertical section $(E_\alpha)_x \subseteq Y$ contains a  dense open subset of $Y$.
 There exists $y\in \bigcap \{(E_\alpha)_x : \alpha < \lambda \}$, since  $\lambda < \cov (\Mee_Y)$. So, $(x,y)\in \bigcap \{E_\alpha : \alpha < \lambda \}$.
                 \end{proof} 
  Cases when the product topology  of two Baire space is a Baire space are generalized onto  $\Cee$-cross topologies by Theorem \ref{BBB}.
   This theorem  does not work whenever $\pi(Y)  \geq \cov (\Mee_X)$. Let $D^\lambda$ be a Cantor cube. Obviously, $\pi(D^\lambda)= \lambda$. In \cite{kp} it was  explicitly observed,  compare also D. Fremlin, T. Natkaniec and I. Rec\l aw \cite{fnr},  that:  
  
 \textit{If a set $E\subseteq X\times D^\lambda$ is open and dense with respect to the  product topology, then there exists a meager subset $P\subseteq X$ such that   any vertical section $E_x= \{ y\in D^\lambda: (x,y) \in E \}$ is dense in $D^\lambda$, for each  $x \in X \setminus P$}.
  
   This implies that  the product topology on $ X\times D^\lambda$ satisfies the Baire category theorem, whenever  $X$ is a Baire space. Our results suggest  a list of questions.   For example:
           \begin{problem}
Does the Baire category theorem hold for the cross topology on $X\times D^\lambda$, whenever  $\lambda \geq \cov (\Mee_X)\geq \omega_1$?
              \end{problem}
              \section{The plane equipped with a  $\mathcal C$-cross topology}
From now on we consider the plane $\Ree \times \Ree$ with various $\mathcal C$-cross topologies. The cross topology on  the plane  is not 
    quasiregular. Indeed, every graph of an one-to-one function is closed and nowhere dense with respect to the cross topology, see \cite{hk} Proposition 1.2. There are many 
    one-to-one functions with $\tau$-dense graphs. Any complement of  a such  graph witnesses that the cross topology is not quasiregular, since    $\emptyset \not= V \in \gamma$ implies $ \int_\tau \cl_\gamma V \not= \emptyset$ by Lemma 3.4  in \cite{hw}.
      From this lemma it follows that topologies $\tau$ and $\sigma$ are $\Pi$-related. 
 Therefore $\nwd_\tau = \nwd_\sigma$ and  $\Mee_\tau = \Mee_\sigma$.  However, these equalities do not hold for $\tau$ and  $\gamma$. There hold $\nwd_\tau\subset\nwd_\gamma$ and $\Mee_\tau\subset \Mee_\gamma$. Indeed, suppose that $F=\cl_\tau F\in
\nwd_\tau$. Then $\cl_\gamma F=F$, since $\gamma$ is  finer then $\tau$. Hence $\int_\tau F=\emptyset $. Consequently $F\in\nwd_\gamma$, since $\int_\tau \cl_\gamma F = \emptyset$  Analogically, one verifies that $\Mee_\tau\subset \Mee_\gamma$.  Because  any graph of an one-to-one function belongs to $\nwd_\gamma$, it should be $\nwd_\tau\not=\nwd_\gamma$ and $\Mee_\tau\not= \Mee_\gamma$. 
             \begin{theorem}\label{HHH} If $F \in \nwd_\tau$, then $F$ is nowhere dense with respect to a $\mathcal C$-cross topology $\mu$, whenever the family $\Cee$ contains a $\pi$-base for the product topology $\tau$. 
              \end{theorem}
             \begin{proof}  Let $F \in \nwd_\tau$ and let $\mathcal U \subseteq \mathcal C$ be a $\pi$-base for $\tau$. The union 
$W=\bigcup\{V\in \Uee :V\cap F=\emptyset \}$
is $\tau$-open and $\tau$-dense. It is also $\mu$-open, since it is the union of $\mu$-open sets which belong to $\Uee$. Suppose that $(p,q)\in H=\int_\mu(\Ree\times\Ree\setminus W)$. Then for any $x\in \int_{\Ree}H^q$ the section $W_x$ is not dense in $\Ree$. We have a contradiction with the Kuratowski-Ulam theorem which one  applies with  $W$ and $\tau$. 
             \end{proof}
              \begin{corollary}\label{FFF} If $F \in \mathcal M_\tau$, then $F$ is meager  with respect to a $\mathcal C$-cross topology, whenever the family $\Cee$ contains a $\pi$-base for the product topology $\tau$. \qed
             \end{corollary} 
The next theorem is formulated for the cross topology $\gamma$. However, it holds for any $\Cee$-cross topology which satisfies Theorem \ref{EEE}  and such that graphs of one-to-one functions are nowhere dense.   For Hausdorff spaces $X$ and $Y$  graphs of one-to-one functions are nowhere dense with respect to the cross topology,
see \cite{hk} Proposition 1.2.4.
              \begin{theorem}\label{GGG}
$\cof (\nwd_\gamma)> 2^\omega$ and  $\cof (\mathcal M_\gamma)> 2^\omega$.
              \end{theorem}
              \begin{proof}
Consider a transfinite family  $\{ F_\alpha: \alpha< 2^\omega \} \subseteq \nwd_\gamma.$ Choose $(p_0,q_0)\in\Ree\times\Ree\setminus F_0$. Assume that points $$\{(p_\beta,q_\beta)\in \Ree\times\Ree\setminus F_\beta : \beta < \alpha \} $$ which have been already chosen constitute the graph of an one-to-one function. By Theorem \ref{EEE}, there exists  a meager subset $P_\alpha\subset \Ree$ such that the vertical section $(\Ree\times \Ree \setminus F_\alpha )_x\subseteq \Ree$ is open and dense for each $x\in \Ree\setminus P_\alpha$. Choose a point $p\in \Ree \setminus (P_\alpha\cup \{p_\beta: \beta <\alpha\})$ and a point $q\in(\Ree \times \Ree \setminus F_\alpha)_{p}\setminus \{q_\beta:\beta<\alpha\}$. Put $p=p_\alpha$ and $q=q_\alpha$. The set $\{(p_\alpha,q_\alpha): \alpha<2^\omega\}\subset \nwd_\gamma$ is  contained in no $F_\alpha$. Hence $\cof(\nwd_\gamma)>2^\omega.$ The proof that $\cof (\Mee_\gamma)>2^\omega$ is analogical.
               \end{proof}
                \begin{theorem}\label{JJJ} If $\mu$ is  a $\mathcal C$-cross topology such that  the family $\Cee$ contains a $\pi$-base for the product topology $\tau$, then
  $\cov(\mathcal M_\mu) = \cov(\mathcal M_{\Ree}).$ \end{theorem}
                \begin{proof}  
Let $\Fee \subset\nwd_{\Ree}$ be a family of closed subsets which witnesses $\cov(\Mee_\Ree) =\cov(\nwd_\Ree)$. Put $\Hee =\{\Ree\times V:V\in \Fee \}$.  By the definitions  $\Hee \subseteq\nwd_\tau$. By Theorem \ref{HHH}, one infers  $\Hee \subset\nwd_\mu$.  Since  $\bigcup \Hee = \Ree\times \Ree$, then $\cov(\mathcal M_\mu) \leq \cov(\mathcal M_{\Ree})$.

Assume that  $\{F_\alpha:\alpha<\cov(\Mee_\mu)\}\subset\nwd_\mu$ 
is a family of $\mu$-closed set which witnesses $\cov(\Mee_\mu).$ By Theorem \ref{EEE},  there is a meager set $P_\alpha\subseteq \Ree$ such that for any section $(\Ree\times\Ree\setminus F_\alpha)_x\subseteq\Ree$ is   a dense and open  for each index $\alpha$ and  any point $x\in \Ree \setminus P_\alpha$. Suppose that  $\cov(\Mee_\mu)<\cov(\Mee_\Ree)$. Hence, there exist  a point $x\in\Ree\setminus\bigcup \{P_\alpha:\alpha<\cov(\Mee_\mu)\}$ and a point   $y \in \bigcap\{ (\Ree\times\Ree\setminus F_\alpha)_x:\alpha<\cov(\Mee_\mu)\}$. So,   $(x,y)\notin \bigcup \{F_\alpha:\alpha<\cov(\Mee_\mu)\}$, a contradiction. 
                \end{proof}
                 \begin{theorem}\label{KKK} Let $\mu$ be  a $\mathcal C$-cross topology such that  the family $\Cee$ contains a $\pi$-base for the product topology $\tau$. If $X\notin\Mee_\Ree$, then the square $X \times X$ is not meager with respect to $\mu.$ Moreover, $\non(\mathcal M_\mu) = \non(\mathcal M_{\Ree}).$
            \end{theorem}
            \begin{proof}  Take $X\notin\Mee_\Ree$.
 Suppose that, the square $X \times X$ is meager with respect to $\mu.$ This means  $F_1\cup F_2\cup\ldots =X\times X$, where any set  $F_n\in\nwd_\mu$. Apply  Theorem \ref{EEE} with $X\times X$. There exist meager subsets $P_n\subset X$ such that  any section $(X\times X\setminus F_n)_x\subseteq X$ is open and dense, for each $x\in X\setminus P_n$. Hence, for any point  $x\in X\setminus (P_1\cup P_2\cup\ldots )$ the intersection  $\bigcap \{ (X\times X\setminus F_n)_x : n=1,2, \ldots \}$ is not meager, a contradiction. So,  
if a set $X\subseteq\Ree$  witnesses $\non(\Mee_\Ree)$, then $X\times X$ witnesses $\non(\Mee_\mu).$ This follows $\non(\Mee_\mu)\leq\non(\Mee_\Ree).$ 

From Theorem \ref{HHH} one infers that any not  $\mu$-meager set is not $\tau$-meager, too. Hence, $\non(\Mee_\mu)\geq\non(\Mee_\Ree).$
              \end{proof}
             \begin{problem}   Is $\add(\mathcal M_\gamma) \not=\add(\mathcal M_{\Ree})$   consistent  with ZFC?
              \end{problem}
           \section{Miscellanea of cofinality }
If a Hausdorff space $X$ is  dense in itself, then  $\cof(\nwd_X)\geq \omega_1. $ Indeed, let $U_0,U_1, \ldots $ be an infinite sequence of pairwise disjoint, non-empty and  open subsets of $X$. Assume that $F_0,F_1, \ldots $ is a sequence of nowhere dense subsets. For each $n$ choose a point  $x_n \in U_n \setminus F_n $.  A family of all points $x_n$  is   a nowhere dense subset of $X$ and no $F_n$ contains this family. So,  no sequence of nowhere
 dense subsets witnesses $\cof(\nwd_X)$. But for some $T_1$-spaces it can
 be $\cof(\nwd_X)=\omega_0$. For example, whenever  $X$ is countable and all its co-finite subsets  are open, only.

Let $Y$ be  a  separable dense in itself metric space.
    There holds $\add (\nwd_Y)= \omega_0 = \non(\nwd_Y)$, since $Y$ contains a copy of the rationals as a dense subset.  If $Y\notin \Mee_Y$, then $ \cov( \Mee_Y) = \cov(\nwd_Y)\geq \omega_1$.  Many differences occur for  spaces  constructed under additional set-theoretical axioms, compare \cite{bla}. However, the equality  
 $\cof (\Mee_\Ree)=\cof (\nwd_\Ree)$ was already proved by D. Fremlin \cite{f}, and a simpler proof was given by B. Balcar, F. Hern\'andez-Hern\'andez and M. Hru\v s\'ak \cite{bhh}.   
 From  \cite{bhh}, see  Theorem 1.6 and Fact 1.5, it follows that $\cof (\Mee_\Ree)=\cof (\nwd_ Y)$. 
   Little is known about relations   between $\cof (\Mee_X)$ and $\cof (\nwd_X)$, whenever $X$ is an arbitrary  topological space. 
  For   metric  spaces we extract topological properties used   in \cite{bhh}, see the proof of Fact 1.5, to get the following. 
  \begin{lemma}\label{lm2} Let $X$ be a  dense in itself metric space. If  $Q$ is a dense subset of $X$, then for each $F \in \nwd_X$ there exists $G\in \nwd_Q$ such that    $ F \subseteq \cl_X G$.
  \end{lemma}
  \begin{proof}  Let $B_n$ be the union of all balls with the radius $\frac 1 n$ and with centers  in  $F \in \nwd_X$. Choose maximal, with respect to the inclusion, sets $A_n \subset Q \cap (B_n \setminus B_{n+1})$ such that distances between points of $A_n$ are greater than  $\frac 1 n$. Put $A_1\cup A_2 \cup \ldots = G$. 
  \end{proof}
   \begin{corollary}\label{nwd} Let $X$ be a  dense in itself metric space. If  $Q$ is a dense subset of $X$, then  $\cof(\nwd_Q)  =\cof (\nwd_X)$.
  \end{corollary}
  \begin{proof}   If $Q \subseteq X$, then $\cof(\nwd_Q)  \leq \cof (\nwd_X)$, since this inequality holds for any dense in itself topological space $X$. Lemma \ref{lm2} follows the inverse inequality.    
  \end{proof} 
   Consider $\lambda^\omega$ with the   product topology, where  a  cardinal number $\lambda$ is equipped with the discrete topology.  By the theorem of P. \v St\v ep\'anek and P. Vop\v enka \cite{sv}, compare a general version of this theorem \cite{ks}, one obtains $\omega_1=\add (\Mee_{\lambda^\omega}) =  \cov (\Mee_{\lambda^\omega})= \cov (\nwd_{\lambda^\omega}).$ Any metrizable space    has a $\sigma$-discrete base.  Each selector defined on elements of a such base witnesses $\add (\nwd_{\lambda^\omega}) =\omega_0$, since and $\lambda^\omega$ is a metrizable space. Indeed, any such selector is a countable union of discrete subsets and each discrete subset of a dense in itself $T_1$-space is nowhere dense.   
 Every non-empty open subset of $\lambda^\omega$ contains a family  of the cardinality $\lambda$ which  consists of  non-empty, pairwise disjoint and  open subsets. Hence $\cof (\nwd_{\lambda^\omega}) > \lambda$. Each subset of $\lambda^\omega$ of the cardinality less than $\lambda$ has to be nowhere dense, therefore   $\non(\nwd_{\lambda^\omega})= \lambda$. 
        \begin{theorem} \label{4XX} If    $X$is a metric space such that any non-empty open subset of  $X$ has density $\lambda$,  then $\cof(\nwd_{\lambda^\omega}) = \cof(\nwd_{X}).$
         \end{theorem}
\begin{proof} Let $Y(\lambda )$ be the universal $\sigma$-discrete metric space in the class of all $\sigma$-discrete metric spaces of the cardinality less or equal to $\lambda$.  One can define $Y(\lambda )$ similar as the space $Y(S)$ in \cite{pl} p. 41 or as $Q(\tau )$ in \cite{pr} p. 217.  One can check that  $X$ and $\lambda^\omega$ contain dense homeomorphic copies of  $Y( \lambda) $, compare the proof of Theorem 1 in \cite{pl} or Corollary 7.7 in \cite{pr} or Theorem 1 in \cite{MED}. We are done by Corollary \ref{nwd}. 
\end{proof}
If $\mathbb Q\subset\Ree$ is the set of rational numbers, then the square $\mathbb Q\times\mathbb Q$ is  $\gamma$-dense. Therefore one infers that $\non(\nwd_\gamma)=\add(\nwd_\gamma )=\omega_0$.     
Lemma \ref{lm2} or Corollary \ref{nwd} hold  for the topology $\sigma$, since   $\nwd_\sigma =  \nwd_\tau$. It is impossible to use $\gamma$ instead of $\tau$ in these facts. Indeed, $\cof(\nwd_\gamma)>2^\omega$ implies  $\cof (\nwd_{\mathbb Q\times\mathbb Q}) < \cof (\nwd_\gamma)$, whenever $\mathbb Q\times\mathbb Q$  inherits  its topology  from the topology $\gamma$.

\end{document}